\documentclass[submission]{dmtcs}

\newcommand{\AlternateState}	{ R }

\newcommand{\Prob}	{ {\rm P} }

\newcommand{\Measurable}	{ \sigma}

\newcommand{\ProcessAlphabet}	{\mathcal{A}}
\newcommand{\LocalPast}			{ L^{-}} 
\newcommand{\localpast}			{ l^{-}}
\newcommand{\localpastprime}		{ \lambda}

\newcommand{\LocalFuture}		{ {L}^{+}}

\newcommand{\LocalState}		{ S}
\newcommand{\localstate}		{ s}

\newcommand{\PatchPast}			{ P^{-}}

\newcommand{\PatchFuture}		{ P^{+}}

\newcommand{\indep}			{ \rotatebox{90}{$\models$}}
\newcommand{\nindep}	{\not\hspace{-.05in}\indep}

\newcommand{\GraphPoint}{\left\langle v,t\right\rangle}
\newcommand{\AltGraphPoint}{\left\langle u,s\right\rangle}
\newcommand{\GraphPointA}{\left\langle u,t\right\rangle}
\newcommand{\GraphPointB}{\left\langle v,t\right\rangle}
\newcommand{\GraphPointC}{\left\langle v,s\right\rangle}
\newcommand{\GraphPointPlus}{\left\langle v,t+1\right\rangle}
\newcommand{\GraphPointBMinus}{\left\langle u, t-1\right\rangle}

\usepackage[round]{natbib}
\usepackage{graphics}

\newtheorem{definition}{Definition}
\newtheorem{theorem}{Theorem}
\newtheorem{lemma}{Lemma}
\newtheorem{proposition}{Proposition}
\newtheorem{conjecture}{Conjecture}
\newtheorem{corollary}{Corollary}

\author{Cosma Rohilla Shalizi\thanks{Supported by a grant from the James
    S. McDonnell Foundation}}

\title[Prediction of Random Fields on Networks]{Optimal Nonlinear Prediction of
  Random Fields on Networks}

\address{Center for the Study of Complex Systems, 4485 Randall Laboratory,
  University of Michigan, Ann Arbor, MI 48109 USA}

\keywords{Networks, random fields, sufficient statistics, nonlinear prediction,
  information theory, recursive estimation}

\revision{1} \revised{} \received{11 May 2003} \accepted{}

\begin{document}

\maketitle

\begin{abstract}
It is increasingly common to encounter time-varying random fields on networks
(metabolic networks, sensor arrays, distributed computing, etc.).  This paper
considers the problem of optimal, nonlinear prediction of these fields, showing
from an information-theoretic perspective that it is formally identical to the
problem of finding minimal local sufficient statistics.  I derive general
properties of these statistics, show that they can be composed into global
predictors, and explore their recursive estimation properties.  For the special
case of discrete-valued fields, I describe a convergent algorithm to identify
the local predictors from empirical data, with minimal prior information about
the field, and no distributional assumptions.
\end{abstract}

\tableofcontents

\section{Introduction}

Within the field of complex systems, most interest in networks has focused
either on their structural properties
\citep{MEJN-on-network-structure-and-function}, or on the behavior of known
dynamical systems coupled through a network, especially the question of their
synchronization \citep{Pikovsky-Rosenblum-Kurths-synchronization}.  Statistical
work, ably summarized in the book by \citet{Guyon-random-fields}, has largely
(but not exclusively) focused on characterization and inference for {\em
  static} random fields on networks.  In this paper, I consider the problem of
predicting the behavior of a random field, with unknown dynamics, on a network
of fixed structure.  Adequate prediction involves reconstructing those unknown
dynamics, which are in general nonlinear, stochastic, imperfectly measured, and
significantly affected by the coupling between nodes in the network.  Such
systems are of interest in many areas, including biochemistry
\citep{Bower-Bolouri}, sociology \citep{Young-strategy-structure}, neuroscience
\citep{Dayan-Abbott}, decentralized control
\citep{Siljak-decentralized-control,Mutambara-decentralized} and distributed
sensor systems \citep{Varshney-distributed-detection}.  Cellular automata
\citep{Ilachinski-discrete} are a special case of such systems, where the graph
is a regular lattice.

There are two obvious approaches to this problem of predicting network
dynamics, which is essentially a problem of system identification.  One is to
infer a global predictor, treating entire field configurations as measurements
from a time series; this is hugely impractical, if only because on a large
enough network, no configuration ever repeats in a reasonable-sized data
sample.  The other straightforward approach is to construct a distinct
predictor for each node in the network, treating them as so many isolated
time-series.  While more feasible, this misses any effects due to the coupling
between nodes, which is a serious drawback.  In these contexts, we often know
very little about the causal architecture of the systems we are studying, but
one of the things we do know is that the links in the network are causally
relevant.  In fact, it is often {\em precisely} the couplings which interest
us.  However, node-by-node modeling ignores the effects of the couplings, which
then show up as increased forecast uncertainty at best, and systematic biases
at worst.  We cannot guarantee optimal prediction unless these couplings are
taken into account.  (See the end of Section
\ref{sec:minimal-local-sufficient}.)  I will construct a distinct predictor for
each node in the network, but these predictors will explicitly take into
account the couplings between nodes, and use them as elements in their
forecasts.

By adapting tools from information theory, I construct optimal, nonlinear {\em
  local} statistical predictors for random fields on networks; these take the
form of minimal sufficient statistics.  I describe some of the optimality
properties of these predictors, and show that the local predictors can be
composed into global predictors of the field's evolution.  Reconstructing the
dynamics inevitably involves reconstructing the underlying state space, and
raises the question of determining the state from observation, i.e., the
problem of filtering or state estimation.  There is a natural translation from
the local predictors to a filter which estimates the associated states, and is
often able to do so with {\em no} error.  I establish that it is possible to
make this filter recursive without loss of accuracy, and show that the filtered
field has some nice Markov properties.  In the special case of a discrete
field, I give an algorithm for approximating the optimal predictor from
empirical data, and prove its convergence.

Throughout, my presentation will be theoretical and abstract; for the most
part, I will not deal with practical issues of implementing the method. In
particular, I will slight the the important question of how much data is needed
for reliable prediction.  However, the method {\em has} been successfully
applied to fields on regular lattices (see Section \ref{sec:reconstruction}
below).

I will make extensive use of (fairly basic) information theory and properties
of conditional independence relations, and accordingly will assume some
familiarity with conditional measures.  I will not pay attention to
measure-theoretic niceties, and shall assume that the random field is
sufficiently regular that all the conditional measures I invoke actually exist
and are conditional probability distributions.  Readers, for their part, may
assume all functions to be measurable functions.

The next section of the paper establishes the basic setting, notation, and
preliminary results, the latter mostly taken without proof from standard
references on information theory and conditional independence.  Section
\ref{sec:causal-states} constructs the local predictors and establishes their
main properties.  Section \ref{sec:connections} establishes the results about
transitions between states which are related to recursive filtering.  Section
\ref{sec:reconstruction} discusses an algorithm for identifying the states from
empirical data on discrete fields.

\section{Notation and Preliminaries}
\label{sec:preliminaries}

\subsection{The Graph and the Random Field $X(\GraphPoint)$}

Consider a fixed graph, consisting of a set of nodes or vertices $V$, and
undirected edges between the nodes $E$.  An edge is an ordered pair $(v_1,
v_2)$, indicating that it is possible to go directly from $v_1$ to $v_2$; the
set of edges is thus a binary relation on the nodes.  We indicate the presence
of this relation, i.e. of a direct path, by $v_1 E v_2$; since the edges are
undirected, this implies that $v_2 E v_1$.  There is a path of length $k$
between two nodes if $v_1 E^k v_2$.  In addition to the graph, we have a time
coordinate which proceeds in integer ticks: $t \in {\Bbb N}$ or $\in {\Bbb Z}$.
We need a way to refer to the combination of a vertex and a time; I will call
this a {\em point}, and write it using the ordered pair $\GraphPoint$.

\begin{figure}
\begin{center}
\resizebox{3.0in}{!}{\includegraphics{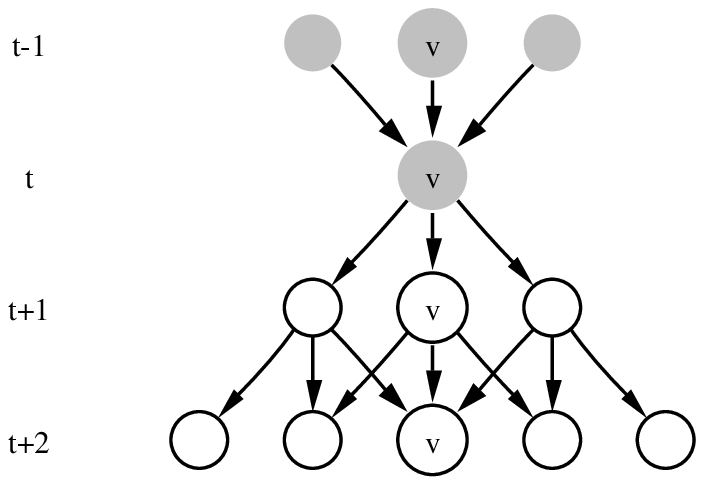}}
\end{center}
\caption{Schematic of the light-cones of node $\GraphPoint$.  Time runs
  vertically downward.  For visual simplicity, vertices are drawn as though the
  graph was a one-dimensional lattice, but no such assumption is necessary.
  The gray circles denote the past light-cone of $\GraphPoint$, and the white
  ones its future light-cone.  Note that $\GraphPoint$ is included in its past
  light-cone, resulting in a slight asymmetry between the two cones.}
\label{st-regions-for-1-cell}
\end{figure}

At each point, we have a random variable $X(\GraphPoint)$, taking values in the
set $\ProcessAlphabet$, a ``standard alphabet'' \citep{Gray-entropy}, such as a
set of discrete symbols or a finite-dimensional Euclidean vector space; this is
the random field. Let $c$ be the maximum speed of propagation of disturbances
in the field.  Now define the {\em past light-cone} of the point $\GraphPoint$
as the set of all points where the field could influence the field at
$\GraphPoint$ (including $\GraphPoint$); likewise the {\em future light-cone}
is all the points whose fields could be influenced by the field at
$\GraphPoint$ (excluding $\GraphPoint$).  We will mostly be concerned with the
configurations in these light-cones, rather than the sets of points themselves;
the configurations in the past and future light-cone are respectively
\begin{eqnarray}
\LocalPast(\GraphPoint) & \equiv & \bigcup_{\tau \geq 0}{\left\{X(\left\langle u,t-\tau\right\rangle)\left|\bigvee_{k=0}^{k=c\tau}{u E^k v}\right.\right\}} ~\mathrm{and}\\
\LocalFuture(\GraphPoint) & \equiv & \bigcup_{\tau \geq
  1}{\left\{X(\left\langle u,t+\tau\right\rangle)\left|\bigvee_{k=0}^{k=c\tau}{v E^k u}\right.\right\}} ~.
\end{eqnarray}
Figure \ref{st-regions-for-1-cell} is a schematic illustration of the past and
future light-cones\footnote{The term ``light-cone'' is used here in analogy
  with relativistic physics, but it's just an analogy.}  There is a certain
distribution over future light-cone configurations, conditional on the
configuration in the past; this is
$\Prob(\LocalFuture(\GraphPoint)|\LocalPast(\GraphPoint) = \localpast)$, which
I shall abbreviate $\Prob(\LocalFuture|\localpast)$.  Note that, in general,
the light-cones of distinct vertices have completely different shapes, and so
neither their configurations nor the distribution over those configurations are
comparable.

\subsection{Mutual Information}

The mutual information between two random variables $X$ and $Y$ is
\begin{eqnarray}
I[X;Y] & \equiv & \mathbf{E}_{X,Y}\left[\log_2{\frac{\Prob(X = x, Y = y)}{\Prob(X = x)\Prob(Y= y)}}\right] ~,
\end{eqnarray}
where $\mathbf{E}$ is mathematical expectation, and $\Prob$ is the
probability mass function for discrete variables, and the probability density
function for continuous ones.  That is, it is the expected logarithm of the
ratio between the actual joint distribution of $X$ and $Y$, and the product of
their marginal distributions.  $I[X;Y] \geq 0$, and $I[X;Y] = 0$ if and only if
$X$ and $Y$ are independent.

The conditional mutual information $I[X;Y|Z]$, is
\begin{eqnarray}
I[X;Y|Z] & \equiv & \mathbf{E}_{X,Y,Z}\left[\log_2{\frac{\Prob(X = x, Y = y|Z=z)}{\Prob(X = x|Z=z)\Prob(Y= y|Z=z)}}\right]
\end{eqnarray}
Conditional mutual information is also non-negative.

We will need only two other properties of mutual information.  The first is
called the ``data-processing inequality''.  For any function $f$,
\begin{eqnarray}
\label{eqn:data-processing}
I[f(X);Y] & \leq & I[X;Y] ~\mathrm{and} \\
\label{eqn:conditional-data-processing}
I[f(X);Y|Z] & \leq & I[X;Y|Z] ~.
\end{eqnarray}
The other is called the ``chain rule''.
\begin{eqnarray}
\label{eqn:chain-rule}
I[X,Z;Y] & = & I[X;Y] + I[Z;Y|X]
\end{eqnarray}

For more on the nature and uses of mutual information, and proofs of Equations
\ref{eqn:data-processing}--\ref{eqn:chain-rule}, see any text on information
theory, e.g.\ \citet{Gray-entropy}.

\subsection{Conditional Independence}

Two random variables, $X$ and $Y$, are conditionally independent given a third,
$T$, in symbols $X \indep Y | T$, when
\begin{eqnarray}
\label{cond-ind-if-joint-prob-factors}
\Prob(X,Y|T) & = & \Prob(X|T)\Prob(Y|T) ~.
\end{eqnarray}
$X \indep Y|T$ if and only if
\begin{eqnarray}
\label{cond-ind-if-other-var-irrelevant}
\Prob(Y|X,T) & = & \Prob(Y|T) ~\mathrm{and}\\
\label{eqn:conditional-independence-and-mutual-information}
I[X;Y|T] & = & 0 ~.
\end{eqnarray}

Conditional independence has many implicative properties; we will use some of
the standard ones, proofs of which may be found in any book on graphical models
\citep{Lauritzen-graphical-models,Pearl-causality,Spirtes-Glymour-Scheines}.
\begin{eqnarray}
\label{intersection}
(A \indep B|CD) \wedge (A \indep D|CB) & \Rightarrow & (A \indep BD|C)\\
\label{weak-union}
(A \indep BC | D) & \Rightarrow & (A \indep B | CD)\\
\label{contraction}
(A \indep B| C) \wedge (A \indep D | CB) & \Rightarrow & (A \indep BD | C)
\end{eqnarray}
Sadly,
\begin{eqnarray}
(A \indep B) & \not\Rightarrow & (A \indep B|C)\\
(A \indep B|C) & \not\Rightarrow & (A \indep B|CD)
\end{eqnarray}

The following property, while not included in most lists of conditional
independence properties, is of some use to us:
\begin{eqnarray}
(A \indep B|C) & \Rightarrow & (A \indep f(B)|C)
\label{cond-ind-of-variable-implies-cond-ind-of-function}
\end{eqnarray}
for any nonrandom function $f$.  It is a direct consequence of the
non-negativity of conditional mutual information, and Equations
\ref{eqn:conditional-data-processing} and
\ref{eqn:conditional-independence-and-mutual-information}.

\subsection{Predictive Sufficient Statistics}

Any function on the past light-cone defines a {\em local statistic}, offering
some summary of the influence of all the points which could affect what happens
at $\GraphPoint$.  A good local statistic conveys something about ``what comes
next'', i.e., $\LocalFuture$.  We can quantify this using information theory.

To be slightly more general, suppose we know the random variable $X$ and wish
to predict $Y$.  Any function $f(X)$ defines a new random variable $F = f(X)$
which is a {\em statistic}.  Because $F$ is a function of $X$, the data
processing inequality (Equation \ref{eqn:data-processing}) says $I[Y;F] \leq
I[Y;X]$.

\begin{definition}[Sufficient Statistic]
\label{defn:sufficiency}
A statistic $F = f(X)$ is {\em sufficient} for $Y$ when $I[Y;F] = I[X;F]$.
\end{definition}
A sufficient statistic, that is to say, is one which is as informative as the
original data.  While I will not elaborate on this here, it is very important
to remember that any prediction method which uses a non-sufficient statistic
can be bettered by one which does.  No matter what the loss function for
prediction, the optimal predictor can always be implemented by an algorithm
which depends solely on a sufficient statistic \citep{Blackwell-Girshick}.

Here are two important criteria, and consequences, of sufficiency.

\begin{proposition}
\label{prop:sufficiency-and-cond-prob}
$F$ is a sufficient statistic, as in Definition \ref{defn:sufficiency}, if and
  only if, $\forall x$,
\begin{eqnarray}
\Prob(Y|X=x) & = & \Prob(Y|F=f(x)) ~.
\end{eqnarray}
\end{proposition}
See \citet[Sections 2.4 and 2.5]{Kullback-info-theory-and-stats}.  (Some
authors prefer to use this as the definition of sufficiency.)

\begin{lemma}
\label{from-cond-ind-to-sufficiency}
$F$ is a sufficient statistic if and only if $X$ and $Y$ are conditionally
independent given $F$, i.e., $X \indep Y | F$.
\end{lemma}

{\em Proof.}  ``Only if'': begin with the observation that
\begin{eqnarray}
\label{eqn:functions-add-no-conditionality}
\Prob(Y|X,f(X)) = \Prob(Y|X) ~,
\end{eqnarray}
no matter what the function\footnote{Let $\mathcal{A} = \Measurable(X)$ be the
  sigma-algebra generated by $X$, and ${\mathcal B} = \Measurable(f(X))$.
  Clearly $\mathcal{B} \subseteq \mathcal{A}$, so $\Measurable(X, f(X)) =
  {\mathcal AB} = {\mathcal A} = \Measurable(X)$.  Thus, the sigma-algebra
  involved when we condition on $X$ is the same as when we condition on $X$ and
  $f(X)$ at once.}.  Hence $\Prob(Y|X,F) = \Prob(Y|X)$.  But if $Y \indep X |
F$, then $\Prob(Y|X,F) = \Prob(Y|F)$.  Hence $\Prob(Y|X) = \Prob(Y|F)$.
``If'': start with $\Prob(Y|F) = \Prob(Y|X)$.  As before, by Equation
\ref{eqn:functions-add-no-conditionality}, since $F = f(X)$, $\Prob(Y|X) =
\Prob(Y|X,F)$.  Hence $\Prob(Y|F) = \Prob(Y|X,F)$, so
(Eq.\ \ref{cond-ind-if-other-var-irrelevant}), ${X \indep Y |F}$.  \qed

While all sufficient statistics have the same predictive power, they are not,
in general, equal in other respects.  In particular, some of them make finer
distinctions among past light-cones than others.  Since these are distinctions
without a difference, we might as well, in the interests of economy, eliminate
them when we find them.  The concept of a {\em minimal} sufficient statistic
captures the idea of eliminating all the distinctions we can get rid off,
without loss of predictive power.

\begin{definition}[Minimal Sufficiency]
\label{definition:minimal-sufficient}
$F$ is a {\em minimal sufficient statistic} for predicting $Y$ from $X$ if and
only if it is sufficient, and it is a function of every other sufficient
statistic.
\end{definition}

\section{Construction and Properties of Optimal Local Predictors}
\label{sec:causal-states}

\subsection{Minimal Local Sufficient Statistics}
\label{sec:minimal-local-sufficient}

We have observed that for each past light-cone configuration $\localpast$ at a
point, there is a certain conditional distribution over future light-cone
configurations, $\Prob(\LocalFuture|\localpast)$.  Let us say that two past configurations, or ``pasts'', are equivalent if they have the same conditional
distribution,
\begin{eqnarray}
\localpast_1 \sim \localpast_2 & \Leftrightarrow & \Prob(\LocalFuture|\localpast_1) = \Prob(\LocalFuture|\localpast_2)
\end{eqnarray}
Let us write the equivalence class of $\localpast$, that is, the set of all
pasts it is equivalent to, as $[\localpast]$.  We now define a local statistic,
which is simply the equivalence class of the past light-cone configuration.

\begin{definition}[Local Causal State]
The local causal state at $\GraphPoint$, written $\LocalState(\GraphPoint)$, is
the set of all past light-cones whose conditional distribution of future
light-cones is the same as that of the past light-cone at $\GraphPoint$.  That
is,
\begin{eqnarray}
\LocalState(\GraphPoint) = \epsilon(\localpast) & \equiv & [\localpast]\\
& = & \left\{{\localpastprime}\left|\Prob(\LocalFuture|\localpastprime) = \Prob(\LocalFuture|\localpast)\right.\right\}
\end{eqnarray}
\end{definition}

The name ``causal state'' was introduced for the analogous construction for
time series by \citet{Inferring-stat-compl}.  I will continue to use this term
here, for two reasons.  First, without getting into the vexed questions
surrounding the nature of causality and the properties of causal relations, it
is clear that the causal states have at least the ``screening'' properties
\citep{Salmon-1984} all authorities on causality require, and may well have the
full set of counter-factual or ``intervention'' properties demanded by
\citet{Pearl-causality} and \citet{Spirtes-Glymour-Scheines}.  (When and
whether they meet the stricter criteria is currently an open question.)
Second, and decisively, these objects need a name, and the previous terms in
the literature, like ``action-test pairs in a predictive state representation''
\citep{predictive-representations-of-state}, ``elements of the statistical
relevance basis'' \citep{Salmon-1984}, or ``states of the prediction process''
\citep{Knight-predictive-view}, are just too awkward to use.

\begin{theorem}[Sufficiency of Local Causal States]
\label{theorem:sufficiency-of-causal-states}
The local causal states are sufficient.
\end{theorem}
{\em Proof}: First, we will find the distribution of future light-cone
configurations conditional on the causal state.  Clearly, this is the average
over the distributions conditional on the past light-cones contained in the
state.  (See \citet[Section 25.2]{Loeve-probability} for details.)  Thus,
\begin{eqnarray}
\Prob(\LocalFuture|\LocalState = \localstate) & = & \int_{\localpastprime \in \epsilon^{-1}(\localstate)}{\Prob(\LocalFuture|\LocalPast = \localpastprime) \Prob(\LocalPast = \localpastprime|\LocalState = \localstate) d\localpastprime}~.
\end{eqnarray}
By construction, $\Prob(\LocalFuture|\LocalPast = \localpastprime)$ is the same
for all $\localpastprime$ in the domain of integration, so, picking out an
arbitrary representative element $\localpast$, we pull that factor out of the
integral,
\begin{eqnarray}
\Prob(\LocalFuture|\LocalState = \localstate) & = & \Prob(\LocalFuture|\LocalPast = \localpast)\int_{\localpastprime \in \epsilon^{-1}(\localstate)}{\Prob(\LocalPast = \localpastprime|\LocalState = \localstate) d\localpastprime}~.
\end{eqnarray}
But now the integral is clearly 1, so, for all $\localpast$,
\begin{eqnarray}
\Prob(\LocalFuture|\LocalState = \epsilon(\localpast)) & = & \Prob(\LocalFuture|\LocalPast = \localpast) ~.
\end{eqnarray}
And so, by Proposition \ref{prop:sufficiency-and-cond-prob}, $\LocalState$ is
sufficient.  \qed

\begin{corollary}
The past and future light-cones are independent given the local causal state.
\label{local-conditional-independence}
\end{corollary}
{\em Proof}: Follows immediately from the conclusion of the theorem and Lemma
\ref{from-cond-ind-to-sufficiency}.  \qed

\begin{corollary}
Let $K$ be the configuration in any part of the future light-cone.  Then, for
any local statistic $\AlternateState$,
\begin{eqnarray}
\label{eqn:maximal-prediction}
I[K;\AlternateState] & \leq I[K;\LocalState] ~.
\end{eqnarray}
\end{corollary}
This corollary is sometimes useful because $I[\LocalFuture;\AlternateState]$
can be infinite for all non-trivial statistics.

\begin{theorem}[Minimality of Local Causal States]
\label{local-refinement-lemma}
The local causal state is a minimal sufficient statistic.
\end{theorem}
{\em Proof}: We need to show that, for any other sufficient statistic
$\AlternateState = \eta(\LocalPast)$ there is a function $h$ such that
$\LocalState = h(\AlternateState)$.  From Proposition
\ref{prop:sufficiency-and-cond-prob} $\Prob(\LocalFuture|\AlternateState =
\eta(\localpast)) = \Prob(\LocalFuture|\localpast)$.  It follows that
$\eta(\localpast_1) = \eta(\localpast_2)$ only if
$\Prob(\LocalFuture|\localpast_1) = \Prob(\LocalFuture|\localpast_2)$.  This in
turn implies that $\localpast_1 \sim \localpast_2$, and so
$\epsilon(\localpast_1) = \epsilon(\localpast_2)$.  Thus, all histories with a
common value of $\eta(\localpast)$ also have the same value of
$\epsilon(\localpast)$, and one can determine $\LocalState$ from
$\AlternateState$.  Hence the required function exists.  \qed

\begin{corollary}[Uniqueness of the Local Causal States]
\label{cor:uniqueness}
If $\AlternateState = \eta(\LocalPast)$ is a minimal local sufficient
statistic, then $\eta$ and $\epsilon$ define the same partition of the data.
\end{corollary}

{\em Proof}: Since $\LocalState$ is a minimal statistic, it is a function of
$\AlternateState$, and for some function $h$, $\LocalState =
h(\AlternateState)$.  But since $\AlternateState$ is also minimal, there is a
function $g$ such that $\AlternateState = g(\LocalState)$.  Hence
$\epsilon(\localpast_1) = \epsilon(\localpast_2)$, if and only if
$\eta(\localpast_1) = \eta(\localpast_2)$.  \qed

In the introduction, I claimed that if we tried to predict the future of the
network by building a separate model for each vertex, taken in isolation, the
results would generally be sub-optimal.  Theorem \ref{local-refinement-lemma}
and Corollary \ref{cor:uniqueness} vindicate this claim.  The only way
predictions based on the vertex-by-vertex procedure could be as good as ones
based on the full light-cone is if {\em all} of the rest of the light-cone was
{\em always} irrelevant.  This would mean that there was, effectively,
no coupling whatsoever between the vertices.

\subsection{Statistical Complexity of the Field}

Much thought has gone into the problem of defining a measure of complexity that
is not simply a measure of randomness, as are Shannon entropy and the
algorithmic information \citep{Badii-Politi}.  Perhaps the best suggestion is
the one which seems to have originated with \citet{Grassberger-1986}, that the
complexity of a system is the minimal amount of information about its state
needed for optimal prediction.  This suggests, following
\citet{Inferring-stat-compl}, that we identify the complexity of the system
with the amount of information needed to specify its causal state.  Crutchfield
and Young called this quantity the {\em statistical complexity}.

In the case of random fields, it is more appropriate to look at a local,
point-by-point version of this quantity.  That is, the {\em local statistical
  complexity} is
\begin{eqnarray}
C(\GraphPoint) & \equiv & I[\LocalState(\GraphPoint);\LocalPast(\GraphPoint)]
~.
\end{eqnarray}
Dropping the argument, we see that because of Equation
\ref{eqn:data-processing} and the fact that $\LocalState$ is a minimal
statistic,
\begin{eqnarray}
\label{eqn:minimal-cmu}
I[\AlternateState;\LocalPast] & \geq & C
\end{eqnarray}
for any other sufficient statistic $\AlternateState$.  In fact, one can start
with Equations \ref{eqn:maximal-prediction} and \ref{eqn:minimal-cmu}, maximal
predictive power and minimal complexity, as axioms, and from them derive all
the properties of the local causal states \citep{CMPPSS}.

A useful property of the statistical complexity is that it provides an upper
bound on the predictive information.
\begin{eqnarray}
I[\LocalFuture;\LocalPast] & \leq & I[\LocalPast;\LocalState]
\end{eqnarray}
{\em Proof}: Use the chain rule for information, Equation
\ref{eqn:chain-rule}.
\begin{eqnarray}
I[\LocalState,\LocalFuture;\LocalPast] & = & I[\LocalState;\LocalPast] +
I[\LocalFuture;\LocalPast|\LocalState]\\
& = &  I[\LocalState;\LocalPast] ~,
\end{eqnarray}
since Corollary \ref{local-conditional-independence} and Equation
\ref{eqn:conditional-independence-and-mutual-information} imply
$I[\LocalFuture;\LocalPast|\LocalState] = 0$.  Using the chain rule the other
way,
\begin{eqnarray}
I[\LocalState,\LocalFuture;\LocalPast] & = & I[\LocalFuture;\LocalPast] + I[\LocalState;\LocalPast|\LocalFuture]\\
& \geq & I[\LocalFuture;\LocalPast] ~,
\end{eqnarray}
since $I[\LocalState;\LocalPast|\LocalFuture] \geq 0$.

\citet{CMPPSS} give detailed arguments for why $C$ is the right way to measure
complexity; here I will just mention two more recent applications.  First, when
our random field can be interpreted as a macroscopic variable which is a
coarse-graining of an underlying microscopic system, as in statistical
mechanics, then $C$ is the amount of information the macro-variable contains
about the micro-state \citep{What-is-a-macrostate}.  Second, the rate of change
of $C$ over time provides a quantitative measurement of self-organization in
cellular automata \citep{Quant-self-org-in-FN03}.

\subsection{Composition of Global Predictors from Local Predictors}
\label{section:composing-global-from-local}

We have just seen that, if we are interested in the future of an individual
point, we can do optimal prediction with only a knowledge of its local causal
state.  One might worry, however, that in compressing the past light-cone down
to the causal state, we have thrown away information that would be valuable on
a larger scale, that would help us if we wanted to predict the behavior of
multiple vertices.  In the limit, if we wanted to predict the behavior of the
entire network, and had only the local causal states available to us, how badly
would we be hampered?

To address these issues, let us turn our thoughts from a single point to a
connected set of vertices taken at a common time $t$, or a {\em patch}.  The
patch has a past and future light-cone ($\PatchPast$ and $\PatchFuture$,
respectively), which are just the unions of the cones of its constituent
points.  Just as we did for points, we can construct a causal state for the
patch, which we'll call the {\em patch causal state}, which has all the
properties of the local causal states.  We now ask, what is the relationship
between the local causal states of the points in the patch, and the patch
causal state?  If we try to predict the future of the patch using just the
local states, are our predictions necessarily impaired in any way?

The answer, it turns out, is {\em no.}

\begin{lemma}[Patch Composition]
The causal state of a patch at one time is uniquely determined by the
composition of all the local causal states within the patch at that time.
\label{composition-of-patch-states}
\end{lemma}

{\em Proof}: We will show that the composition of local causal states within
the patch is a sufficient statistic {\em of the patch}, and then apply
minimality.

\begin{figure}[t]
\begin{center}
\resizebox{3.0in}{!}{\includegraphics{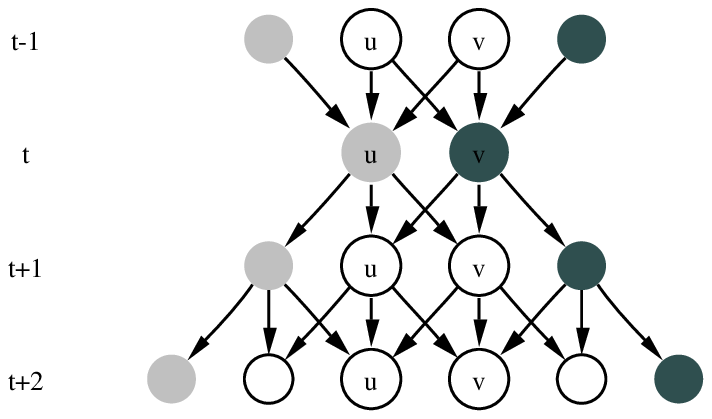}}
\end{center}
\caption{The space-time regions for a two-point patch.  Points which belong
  exclusively to the light-cones of the point on the left ($\GraphPointA$) are
  shaded light grey; those which belong exclusively to the light-cones of the
  other point ($\GraphPointB$) are shaded dark grey.  The areas of overlap
  ($C^{-}$ and $C^{+}$) is white.  Note that, by the definition of
  light-cones, the configuration in $U^{-}$ can have no effect on that in
  $V^{+}$ or vice versa.}
\label{st-regions-for-2-cell-patch}
\end{figure}

\begin{figure}[t]
\begin{center}
\resizebox{!}{1.5in}{\includegraphics{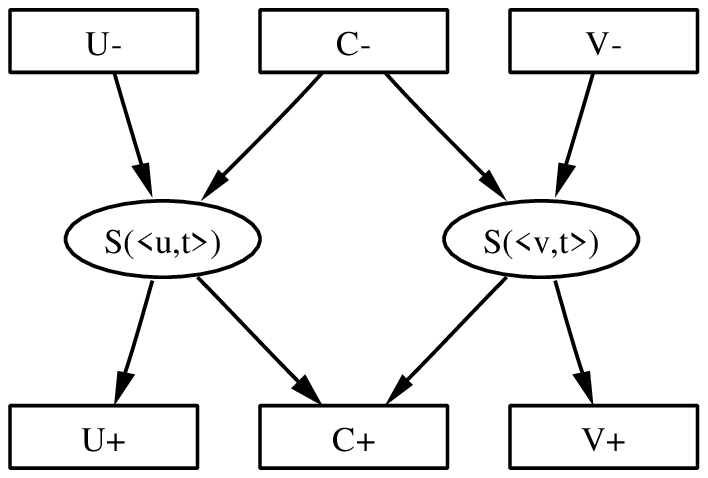}}
\end{center}
\caption{Diagram of effects for the two-node patch.  Arrows indicate the
  direction of influence; causes to effects; the absence of an arrow between
  two nodes indicates an absence of direct influence.}
\label{2-cell-patch-effects-graph}
\end{figure}

Consider first a patch consisting of two spatially adjacent points,
$\GraphPointA$ and $\GraphPointB$.  Define the following variables:
\begin{eqnarray*}
C^{-} & = & \LocalPast(\GraphPointA) \cap \LocalPast(\GraphPointB) \\
U^{-} &  = & \LocalPast(\GraphPointA) \setminus C^{-} \\
V^{-} & = & \LocalPast(\GraphPointB) \setminus C^{-}
\end{eqnarray*}
Thus $\LocalPast(\GraphPointA) = U^{-} \cup C^{-}$, and likewise for
$\LocalPast(\GraphPointB)$.  Define $U^{+}$, $C^{+}$ and $V^{+}$ similarly.
(Figure \ref{st-regions-for-2-cell-patch} gives a picture of these regions.)
Now consider the configurations in these regions.  We may draw a diagram of
effects or influence, Figure \ref{2-cell-patch-effects-graph}, which should not
be confused with the graph of the network.

Corollary \ref{local-conditional-independence} tells us that every path from
$U^{-}$ or $C^{-}$ to $U^{+}$ must go through $\LocalState(\GraphPointA)$.  By
the very definition of light-cones, there cannot be a path linking $V^{-}$ to
$U^{+}$.  Therefore there cannot be a link from $\LocalState(\GraphPointB)$ to
$U^{+}$.  (Such a link would in any case indicate that $U^{+}$ had a dependence
on $C^{-}$ which was not mediated by $\LocalState(\GraphPointA)$, which is
false.)  All of this is true, {\em mutatis mutandis}, for $V^{+}$ as well.

Now notice that every path from variables in the top row --- the variables
which collectively constitute $\PatchPast$ --- to the variables in the bottom
row --- which collectively are $\PatchFuture$ --- must pass through either
$\LocalState(\GraphPointA)$ or $\LocalState(\GraphPointB)$.  The set $Z =
\left\{\LocalState(\GraphPointA), \LocalState(\GraphPointB)\right\}$ thus
``blocks'' those paths.  In the terminology of graphical models
\citep[p.\ 18]{Pearl-causality}, $Z$ {\em d-separates} $\PatchPast$ and
$\PatchFuture$.  But d-separation implies conditional independence ({\em
  ibid.}).  Thus $\PatchPast$ and $\PatchFuture$ are independent given the
composition of $\LocalState(\GraphPointA)$ and $\LocalState(\GraphPointB)$.
But that combination is a function of $\PatchPast$, so Lemma
\ref{from-cond-ind-to-sufficiency} applies, telling us that the composition of
local states is sufficient.  Then Theorem \ref{local-refinement-lemma} tells us
that there is a function from the composition of local states to the patch
causal state.

Now, the reader may verify that this argument would work if one of the two
``nodes'' above was really itself a patch.  That is, if we break a patch down
into a single node and a sub-patch, and we know their causal states, the causal
state of the larger patch is fixed.  Hence, by mathematical induction, if we
know all the local causal states of the nodes within a patch, we have fixed the
patch causal state uniquely.  \qed

\begin{theorem}[Global Prediction]
\label{theorem:global-prediction}
The future of the entire network can be optimally predicted from a knowledge of
all the local causal states at one time.
\end{theorem}
{\em Proof}: Sufficiency, as we've said, implies optimal prediction.  Let us
check that the combination of all the local causal states is a sufficient
statistic for the global configuration.  Apply Lemma
\ref{composition-of-patch-states} to the patch of the entire lattice.  The
proof of the lemma goes through, because it in no way depends on the size or
the shape of the past, or even on the patch being finite in extent.  Since the
patch causal state for this patch is identical with the global causal state, it
follows that the latter is uniquely fixed by the composition of the local
causal states at all points on the lattice.  Since the global causal state is a
sufficient statistic, the combination of all local causal states is too.  \qed

Thus, remarkably enough, the local optimal predictors contain all the
information necessary for global optimal prediction as well.  There doesn't
seem to be any way to make this proof work if the local predictors are not
based on the full light-cones.

\section{Structure of the Field of States $\LocalState(\GraphPoint)$}
\label{sec:connections}

Let us take a moment to recap.  We began with a dynamic random field $X$ on the
network.  From it, we have constructed another random field on the the same
network, the field of the local causal states $S$.  If we want to predict $X$,
whether locally (Theorem \ref{theorem:sufficiency-of-causal-states}) or
globally (Theorem \ref{theorem:global-prediction}), it is sufficient to know
$S$.  This situation resembles attractor reconstruction in nonlinear dynamics
\citep{Kantz-Schreiber}, state-space models in time series analysis
\citep{Durbin-Koopman-state-space-methods}, and hidden Markov models in signal
processing \citep{Elliot-et-al-HMM}.  In each case, it helps, in analyzing and
predicting our observations, to regard them as distorted measurements of
another, unseen set of state variables, which have their own dynamics.  Let us,
therefore, call all these things ``hidden-state models''.

The hidden state space always has a more tractable structure than the
observations, e.g., the former is Markovian, or even deterministic, but the
latter is not.  Now, usually the nice structure is simply {\em demanded} by us
a priori, and it is an empirical question whether the demand can be met,
whether any hidden-state model with that structure can adequately account for
our observations.  However, in the case of time series, one can {\em construct}
hidden states, analogous to those built in Section \ref{sec:causal-states}, and
show that these {\em always} have nice structural properties: they are always
homogeneous Markov processes, and so forth
\citep{Knight-predictive-view,CMPPSS}.  This analogy gives us reason to hope
that the causal states we have constructed always have nice properties, which
is indeed the case.

In this section, I am going to establish some of the structural properties of
the field of causal states, which involve the relations between states at
different points.  Section \ref{sec:transitions} shows how the state at one
point may be used to partially determine the state at other points.  Section
\ref{sec:recursion} applies those results to the problem of designing a filter
to estimate the state field $S$ from the observation field $X$.  Finally,
Section \ref{sec:Markov} considers the Markov properties of the $S$ field.
Note that, if we tried to build a hidden-state model for each vertex
separately, the result would lack these {\em spatial} properties, as well as
being a sub-optimal predictor (see the end of Section
\ref{sec:minimal-local-sufficient}).

\subsection{Transition Properties}
\label{sec:transitions}

The basic idea motivating this section is that the local causal state at
$\GraphPoint$ should be an adequate summary of $\LocalPast(\GraphPoint)$ for
{\em all} purposes.  We have seen that this is true for both local and
non-local prediction.  Here we consider the problem of determining the state at
another point $\AltGraphPoint$, $s \geq t$.  In general, the past light-cones
of $\GraphPoint$ and $\AltGraphPoint$ will overlap.  If we have
$\LocalState(\GraphPoint)$, and want $\LocalState(\AltGraphPoint)$, do
we need to know the actual contents of the overlapping region, or can we get
away with just knowing $\LocalState(\GraphPoint)$?

It seems clear (and we will see that it is true) that determining
$\LocalState(\AltGraphPoint)$ will require knowledge of the ``new''
observations, the ones which are in the past light-cone of $\AltGraphPoint$ but
not in $\LocalPast(\GraphPoint)$.  This data is relevant to $\AltGraphPoint$,
but inaccessible at $\GraphPoint$, and couldn't be summarized in
$\LocalState(\GraphPoint)$.  This data also needs a name; let us call it the
{\em fringe} seen when moving from $\GraphPoint$ to $\AltGraphPoint$.  (See
Figures \ref{st-regions-for-2-cell-patch} and
\ref{st-regions-for-time-forward-transition}.)

The goal of this section is to establish that $\LocalState(\AltGraphPoint)$ is
completely determined by $\LocalState(\GraphPoint)$ and the fringe.  The way I
do this is to show that $\LocalFuture(\AltGraphPoint)$ and
$\LocalPast(\AltGraphPoint)$ are independent, conditional on
$\LocalState(\GraphPoint)$ and the fringe.  This means that the latter two,
taken together, are a sufficient statistic, and then I invoke Theorem
\ref{local-refinement-lemma}.  An alternate strategy would be to consider a
transducer whose inputs are fringes, read as the transducer moves across the
graph, and whose outputs are local states.  Our problem then would be to show
that the transducer's transition rules can always be designed to ensure that
the states returned track the actual states.  This way of framing the problem
opens up valuable connections to the theory of spatial automata and spatial
languages \citep{Two-D-Patterns}, but I prefer a more direct approach by way of
conditional independence.

Even so, the proof is frankly tedious.  First, I show that the desired property
holds if we consider successive states at the same vertex.  Next, I show that
it holds for simultaneous states at neighboring vertices.  Then I extend those
results to relate states at points connected by an arbitrary path in space and
time.  Finally, I show that the result is independent of the precise path
chosen to link the points.  Sadly, I have not been able to find a simpler
argument.

\subsubsection{Temporal Transitions}

\begin{figure}[t]
\begin{center}
\resizebox{3.0in}{!}{\includegraphics{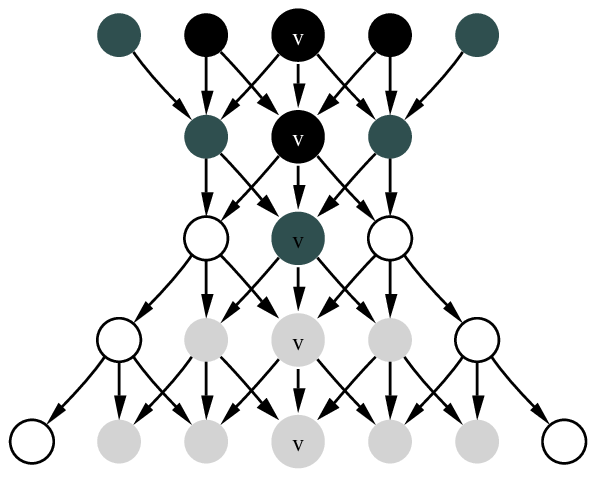}}
\end{center}
\caption{Space-time regions for the time-forward transition from $\GraphPoint$
  to $\GraphPointPlus$.  Black points: $\LocalPast(\GraphPoint)$, the past
  light-cone of $\GraphPoint$.  Dark grey: $N^{-}$, the part of
  $\LocalPast(\GraphPointPlus)$, the past light-cone of $\GraphPointPlus$,
  outside the past light-cone of $\GraphPoint$.  Light grey:
  $\LocalFuture(\GraphPointPlus)$, the future light-cone of $\GraphPointPlus$.
  White: $M^{+}$, the part of $\LocalFuture(\GraphPoint)$ outside
  $\LocalFuture(\GraphPointPlus)$.}
\label{st-regions-for-time-forward-transition}
\end{figure}

We want to move forward in time one step, while staying at the same vertex.
Call the point we start at $\GraphPoint$, and its successor $\GraphPointPlus$.
The whole of the new future light cone is contained inside the old future light
cone, and vice versa for the past cones.  So let's define the following
variables:
\begin{eqnarray*}
N^{-} & = & \LocalPast(\GraphPointPlus) \setminus \LocalPast(\GraphPoint) \\
M^{+} & = & \LocalFuture(\GraphPoint) \setminus \LocalFuture(\GraphPointPlus) ~;
\end{eqnarray*}
$N^{-}$ is the fringe.  (Figure
\ref{st-regions-for-time-forward-transition} pictures these regions.)

\begin{lemma}
The local causal state at $\GraphPointPlus$ is a function of the local causal
state at $\GraphPoint$ and the time-forward fringe $N^{-}$.
\label{lemma:temporal-determinism}
\end{lemma}

\begin{figure}[t]
\begin{center}
\resizebox{!}{2.0in}{\includegraphics{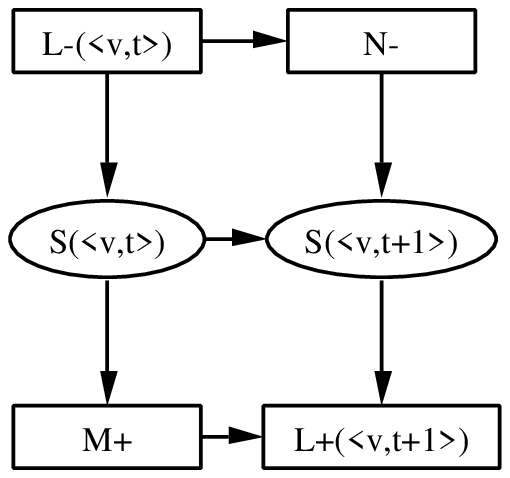}}
\end{center}
\caption{Influence diagram for the variables involved in a time-forward
  transition.}
\label{time-forward-effects-graph}
\end{figure}

{\em Proof.}  Start by drawing the diagram of effects (Figure
\ref{time-forward-effects-graph}).

$M^{+}$ and $\LocalFuture(\GraphPointPlus)$ jointly constitute
$\LocalFuture(\GraphPoint)$, so there must be paths from
$\LocalState(\GraphPoint)$ to both of them.  Now,
$\LocalState(\GraphPointPlus)$ renders $\LocalPast(\GraphPointPlus)$ and
$\LocalFuture(\GraphPointPlus)$ conditionally independent.  Hence it should
d-separate them in the graph of effects.  But $\LocalPast(\GraphPoint)$ is part
of $\LocalPast(\GraphPointPlus)$ and has a direct path to
$\LocalState(\GraphPoint)$.  This means that there cannot be a direct path from
$\LocalState(\GraphPoint)$ to $\LocalFuture(\GraphPointPlus)$; rather, the path
must go through $\LocalState(\GraphPointPlus)$.  (We indicate this in the graph
by a dotted arrow from $\LocalState(\GraphPoint)$ to
$\LocalFuture(\GraphPointPlus)$.  Similarly, $\LocalPast(\GraphPoint)$
certainly helps determine $\LocalState(\GraphPointPlus)$, but it need not do so
directly.  In fact, it cannot: $\LocalState(\GraphPoint)$ must d-separate
$\LocalPast(\GraphPoint)$ and $\LocalFuture(\GraphPoint)$, i.e., must
d-separate $\LocalPast(\GraphPoint)$ from $\LocalFuture(\GraphPointPlus)$ and
$M^{+}$.  Hence the influence of $\LocalPast(\GraphPoint)$ on
$\LocalState(\GraphPointPlus)$ must run through $\LocalState(\GraphPoint)$.
(We indicate this, too, by a dotted arrow from $\LocalPast(\GraphPoint)$ to
$\LocalState(\GraphPointPlus)$.)

Now it is clear that the combination of $\LocalState(\GraphPoint)$ and
$N^{-}$ d-separates $\LocalPast(\GraphPointPlus)$ from
$\LocalFuture(\GraphPointPlus)$, and hence makes them conditionally
independent.  But now the usual combination of Lemma
\ref{from-cond-ind-to-sufficiency} and Theorem \ref{local-refinement-lemma}
tell us that there's a function from $\LocalState(\GraphPoint), N^{-}$
to $\LocalState(\GraphPointPlus)$.  \qed

\subsubsection{Spatial Transitions}

\begin{lemma}
Let $\GraphPointA$ and $\GraphPointB$ be simultaneous, neighboring
points.  Then $\LocalState(\GraphPointB)$ is a function of
$\LocalState(\GraphPointA)$ and the fringe in the direction from
$\GraphPointA$ to $\GraphPointB$, $V^{-}$.
\label{lemma:spatial-determinism}
\end{lemma}

Here the breakdown of the past and future light-cone regions is the same as
when we saw how to compose patch causal states out of local causal states in
Section \ref{section:composing-global-from-local}, as is the influence diagram;
we'll use the corresponding terminology, too.  (See Figures
\ref{st-regions-for-2-cell-patch} and \ref{2-cell-patch-effects-graph},
respectively.)  What we hope to show here is that conditioning on the
combination of $\LocalState(\GraphPointA)$ and $V^{-}$ makes
$\LocalFuture(\GraphPointB)$ independent of $V^{-}$ and $C^{-}$.
Unfortunately, as the reader may verify by inspecting the diagram, our
conditional variables no longer d-separate the other variables (since they have
an unblocked connection through $\LocalState(\GraphPointB)$).  All is not lost,
however: d-separation implies conditional independence, but not conversely.

Abbreviate the pair of variables $\left\{\LocalState(\GraphPointA),
V^{-}\right\}$ by $Z$.  Now, $\LocalState(\GraphPointB)$ is a (deterministic)
function of $C^{-}$ and $V^{-}$.  Hence it is also a function of $Z$ and
$C^{-}$.  Thus $\Prob(V^{+}|\LocalState(\GraphPointB), Z, C^{-}) =
\Prob(V^{+}|Z, C^{-})$.  But this tells us that
\begin{equation}
V^{+} \indep \LocalState(\GraphPointB) | Z, C^{-} 
\end{equation}
From d-separation, we also have
\begin{equation}
V^{+} \indep C^{-} | Z, \LocalState(\GraphPointB)
\end{equation}
Applying Eq.\ \ref{intersection},
\begin{equation}
V^{+} \indep \LocalState(\GraphPointB), C^{-} | Z
\end{equation}
Applying Eq.\ \ref{weak-union},
\begin{equation}
V^{+} \indep C^{-} | Z
\end{equation}
Since $Z = Z, V^{-}$,
\begin{equation}
V^{+} \indep C^{-} | Z, V^{-} 
\end{equation}
The following conditional independence is odd-looking, but trivially true:
\begin{equation}
V^{+} \indep V^{-} | Z
\end{equation}
And it, along with Eq.\ \ref{contraction}, gives us
\begin{equation}
V^{+} \indep C^{-}, V^{-} | Z
\end{equation}

A similar train of reasoning holds for $C^{+}$.  Thus, the entire
future light cone of $\GraphPointB$ is independent of that point's past light
cone, given $\LocalState(\GraphPointA)$ and $V^{-}$.  This tells us that
$\left\{\LocalState(\GraphPointA), V^{-}\right\}$ is sufficient for
$\LocalFuture(\GraphPointB)$, hence $\LocalState(\GraphPointB)$ is a function
of it.  \qed

\subsubsection{Arbitrary Transitions}

\begin{lemma}
Let $\GraphPointA$ and $\GraphPointC$ be two points, $s \geq t$.  Let $\Gamma$
be a spatio-temporal path connecting the two points, arbitrary except that it
can never go backwards in time.  Let $F_{\Gamma}$ be the succession of fringes
encountered along $\Gamma$.  Then $\LocalState(\GraphPointC)$ is a function of
$\LocalState(\GraphPointA)$, $\Gamma$ and $F_{\Gamma}$,
\begin{eqnarray}
\LocalState(\GraphPointC) & = & g(\LocalState(\GraphPointA), \Gamma, F_{\Gamma})
\end{eqnarray}
for some function $g$.
\label{lemma:path-determinism}
\end{lemma}
{\em Proof.}  Apply Lemma \ref{lemma:temporal-determinism} or
\ref{lemma:spatial-determinism} at each step of $\Gamma$.  \qed

\begin{theorem}
Let $\GraphPointA$ and $\GraphPointC$ be two points as in the previous
lemma, and let $\Gamma_1$, $\Gamma_2$ be two paths connecting them, and
$F_{\Gamma_1}$ and $F_{\Gamma_2}$ their fringes, all as in the previous lemma.
Then the state at $\GraphPointC$ is independent of which path was taken to
reach it,
\begin{eqnarray}
g(\LocalState(\GraphPointA), \Gamma_1, F_{{\Gamma}_1}) & = & g(\LocalState(\GraphPointA), \Gamma_2, F_{{\Gamma}_2})~.
\end{eqnarray}
\end{theorem}
{\em Proof.}  Suppose otherwise.  Then either the state we get by going along
$\Gamma_1$ is wrong, i.e., isn't $\LocalState(\GraphPointC)$, or the state
we get by going along $\Gamma_2$ is wrong, or both are.
\begin{eqnarray}
\LocalState(\GraphPointC) \neq g(\LocalState(\GraphPointA), \Gamma_1, F_{{\Gamma}_1}) & \vee & \LocalState(\GraphPointC) \neq g(\LocalState(\GraphPointA), \Gamma_2, F_{{\Gamma}_2}) \\
\LocalFuture(\GraphPointC) \nindep \LocalPast(\GraphPointC) | \LocalState(\GraphPointA), \Gamma_1, F_{{\Gamma}_1} & \vee & \LocalFuture(\GraphPointC) \nindep \LocalPast(\GraphPointC) |\LocalState(\GraphPointA), \Gamma_2, F_{{\Gamma}_2} \\
\neg (\LocalFuture(\GraphPointC) \indep \LocalPast(\GraphPointC) | \LocalState(\GraphPointA), \Gamma_1, F_{{\Gamma}_1} & \wedge & \LocalFuture(\GraphPointC) \indep \LocalPast(\GraphPointC) |\LocalState(\GraphPointA), \Gamma_2, F_{{\Gamma}_2})
\end{eqnarray}
But, by Lemma \ref{lemma:path-determinism}, ${\LocalFuture(\GraphPointC) \indep
  \LocalPast(\GraphPointC) | {\LocalState(\GraphPointA), \Gamma_1,
    F_{{\Gamma}_1}}}$ and ${\LocalFuture(\GraphPointB) \indep
  \LocalPast(\GraphPointC) |{\LocalState(\GraphPointA), \Gamma_2,
    F_{{\Gamma}_2}}}$.  Hence transitions must be path-independent.  \qed

\subsection{Recursive Filtering}
\label{sec:recursion}

Because the causal states are logical constructions out of observational data,
in principle the state at any point can be determined exactly from looking at
the point's past light-cone.  It may be, however, that for some fields one
needs to look infinitely far back into the past to fix the state, or at least
further back than the available data reaches.  In such cases, one will
generally have not exact knowledge of the state, but either a set-valued type
estimate (i.e., $\LocalState \in \mathcal{S}$, for some set of states
$\mathcal{S}$), or a distribution over states, supposing you have a meaningful
prior distribution.

A naive state-estimation scheme, under these circumstances, would produce an
estimate for each point separately.  The transition properties we have just
proved, however, put {\em deterministic} constraints on which states are
allowed to co-exist at different points.  In particular, we can narrow our
estimates of the state at each point by requiring that it be consistent with
our estimates at neighboring points.  Applied iteratively, this will lead to
the tightest estimates which we can extract from our data.  If we can fix the
state at even one point exactly, then this will propagate to all points at that
time or later.

More generally, rather than first doing a naive estimate and then tightening
it, we can incorporate the deterministic constraints directly into a recursive
filter.  Under quite general circumstances, as time goes on, the probability
that such a filter will {\em not} have fixed the state of at least one point
goes to zero, and once it fixes on a state, it stays fixed.  Such filters can
be implemented, at least for discrete fields, by means of finite-state
transducers; see \citet[Chapter 10]{CRS-thesis} for examples (further
specialized to discrete fields on regular lattices).

\subsection{Markov Properties}
\label{sec:Markov}

Recursive estimation is a kind of Markov property \citep{Nevelson-Hasminskii},
and the fact that it can work exactly here is very suggestive.  So, too, is the
fact that the past and the future light-cones are independent, conditional on
the causal state.  It would be very nice if the causal states formed a Markov
random field; for one thing, we could then exploit the well-known machinery for
such fields.  For instance, we would know, from the equivalence between Gibbs
distributions and Markov fields \citep{Guyon-random-fields}, that there was an
effective potential for the interactions between the states across the network.

\begin{definition}[Parents of a Local Causal State]
The parents of the local causal state at $\GraphPoint$ are the causal states at
all points which are one time-step before $\GraphPoint$ and inside its past
light-cone:
\begin{eqnarray}
A(\GraphPoint) & \equiv & \left\{\LocalState(\GraphPointBMinus) \left| (u=v) \vee (uEv) \right. \right\}
\end{eqnarray}
\end{definition}

\begin{lemma}
The local causal state at a point, $\LocalState(\GraphPoint)$, is independent
of the configuration in its past light cone, given its parents.
\begin{equation}
\LocalState(\GraphPoint) \indep \LocalPast(\GraphPoint) | A(\GraphPoint)
\end{equation}
\end{lemma}
{\em Proof.}  $\GraphPoint$ is in the intersection of the future light cones of
all the node in the patch at $t-1$.  Hence, by the arguments given in the proof
of the composition theorem, it is affected by the local states of all those
nodes, {\em and by no others}.  In particular, previous values of the
configuration in $\LocalPast(\GraphPoint)$ have no direct effect; any influence
must go through those nodes.  Hence, by d-separation,
$\LocalState(\GraphPoint)$ is independent of $\LocalPast(\GraphPoint)$.  \qed

\begin{theorem}[Temporal Markov Property]
The local causal state at a point, $\LocalState(\GraphPoint)$, is independent
of the local causal states of points in its past light cone, given its parents.
\end{theorem}
{\em Proof.}  By the previous lemma, $\LocalState(\GraphPoint)$ is
conditionally independent of $\LocalPast(\GraphPoint)$ given its parents.  But
the local causal states in its past light cone are a function of
$\LocalPast(\GraphPoint)$.  Hence by Equation
\ref{cond-ind-of-variable-implies-cond-ind-of-function},
$\LocalState(\GraphPoint)$ is also independent of those local states.  \qed

Comforting though that is, we would really like a stronger Markov property,
namely the following.

\begin{conjecture}[Markov Field Property]
The local causal states form a Markov field in space-time.
\end{conjecture}
{\em Argument for why this is plausible.}  We've seen that, temporally
speaking, a Markov property holds: given a patch of nodes at one time, they are
independent of their past light cone, given their causal parents.  What we need
to add for the Markov field property is that, if we condition on {\em present}
neighbors of the patch, as well as the parents of the patch, then we get
independence of the states of all points at time $t$ or earlier.  It's
plausible that the simultaneous neighbors are informative, since they are also
causal descendants of causal parents of the patch.  But if we consider any more
remote node at time $t$, its last common causal ancestor with any node in the
patch must have been before the immediate parents of the patch, and the effects
of any such local causal state are screened off by the parents.

Unfortunately, this is not really rigorous.  In the somewhat specialized case
of discrete fields on regular lattices, a number of systems have been checked,
and in all cases the states have turned out to be Markov random fields.

\section{State Reconstruction for Discrete Fields}
\label{sec:reconstruction}

The local causal states are a particular kind of hidden-state model; they
combine optimal prediction properties (Section \ref{sec:causal-states}) with
the nice structural properties of hidden-state models (Section
\ref{sec:connections}).  It is only human to be tempted to use them on actual
systems.  In the special case of discrete-valued fields, one can reliably
identify the causal states from empirical data, with minimal assumptions.  This
section describes an algorithm for doing so, building on earlier procedures for
cellular automata \citep{Quant-self-org-in-FN03}.  As always, I take structure
of the network to be known and constant.

Assume that for each point we have estimates of the conditional probability
distribution $\Prob(\LocalFuture, \LocalPast)$ over light-cones for each point.
(These could come, for instance, from the empirical distribution in an ensemble
of networks with the same structure, or from the same network observed over
time, if certain ergodicity assumptions can be made.)  The resulting
conditional distributions, $\Prob(\LocalFuture|\localpast)$, can be treated as
multinomials.  We then cluster past configurations, point by point, based on
the similarity of their conditional distributions.  We cannot expect that the
estimated distributions will match exactly, so we employ a standard test, e.g.
$\chi^2$, to see whether the discrepancy between estimated distributions is
significant.  These clusters are then the estimated local causal states.  We
consider each cluster to have a conditional distribution of its own, equal to
the weighted mean of the distributions of the pasts it contains.

As a practical matter, we need to impose a limit on how far back into the past,
or forward into the future, the light-cones are allowed to extend --- their
depth.  Also, clustering cannot be done on the basis of a true equivalence
relation.  Instead, we list the past configurations
$\left\{\localpast_i\right\}$ in some arbitrary order.  We then create a
cluster which contains the first past, $\localpast_1$.  For each later past,
say $\localpast_i$, we go down the list of existing clusters and check whether
$\Prob_t(\LocalFuture|\localpast_i)$ differs significantly from each cluster's
distribution.  If there is no difference, we add $\localpast_i$ to the first
matching cluster and update the latter's distribution.  If $\localpast_i$ does
not match any existing cluster, we make a new cluster for $\localpast_i$.  (See
Figure \ref{algorithm} for pseudo-code.)  As we give this procedure more and
more data, it converges in probability on the correct set of causal states,
independent of the order in which we list past light-cones (see below).  For
finite data, the order of presentation matters, but we finesse this by
randomizing the order.

\begin{figure}[t]
\begin{tabbing}
\texttt{U} $\leftarrow$ list of all pasts in random order\\
Move the first past in \texttt{U} to a new state\\
for \= each \texttt{past} in \texttt{U}\\
\> \texttt{noMatch} $\leftarrow$ TRUE\\
\> \texttt{state} $\leftarrow$ first state on the list of states\\
\> while \= (\texttt{noMatch} and more states to check)\\
\> \> \texttt{noMatch} $\leftarrow$ (Significant difference between \texttt{past} and \texttt{state}?)\\
\> \> if \= (\texttt{noMatch})\\
\> \> \> \texttt{state} $\leftarrow$ next state on the list\\
\> \> else \=\\
\> \> \> Move \texttt{past} from \texttt{U} to \texttt{state}\\
\> \> \> \texttt{noMatch} $\leftarrow$ FALSE\\
\> if \= (\texttt{noMatch})\\
\> \> make a new state and move \texttt{past} into it from \texttt{U}
\end{tabbing}
\caption{\label{algorithm} Algorithm for grouping past light-cones into estimated states}
\end{figure}

Suppose that the past and the future light-cones are sufficiently deep that
they suffice to distinguish the causal states, i.e., that if we had the exact
distribution over light-cones, the true causal states would coincide exactly
with those obtained from the distribution over limited-depth light-cones.
Then, conditioning on the limited-depth past cones makes futures independent of
the more remote past.  Indeed, every time we examine the future of a given
past, we take an {\em independent} sample from an unchanging distribution over
futures.  Thus, the strong law of large numbers tells us that the empirical
probability of any future configuration will converge on the true probability
almost surely.  Since, with finite future depth, there are only finitely many
possible future configurations, the conditional distribution as a whole also
converges almost surely.  And, if there are only finitely many vertices, we get
over-all convergence of all the necessary conditional distributions.

For this analysis to hold, we need to know three things.
\begin{enumerate}
\item The structure of the graph.
\item The constant $c$ which is the maximum speed at which information can
  propagate.
\item The minimum necessary depth of past and future cones.
\end{enumerate}
Note that if we over-estimate (2) and (3), we will not really harm ourselves,
but under-estimates {\em are} harmful, since in general they will miss useful
bits of predictive information.

Note that this analysis, like the algorithm, is somewhat crude, in that it
doesn't make use of any of the properties of the causal states we established
earlier.  In the analogous case of time series, using those properties can
greatly speed up the reconstruction process, and allows us to not merely prove
convergence, but to estimate the rate of convergence \citep{AfPDiTS}.  Probably
something like that could be done here, with the additional complication that
each vertex must get its own set of states.

\section{Conclusion}

The purpose of this paper has been to define, mathematically, optimal nonlinear
predictors for a class of complex systems, namely dynamic random fields on
fixed, undirected graphs.  Starting with the basic idea of maximizing the
predictive information, we constructed the {\em local causal states} of the
field, which are minimal sufficient statistics, and so the simplest possible
basis for optimal prediction.  We then examined how to combine these states for
non-local prediction, and the structure of the causal-state field.  The last
section described an algorithm for reconstructing the states of discrete
fields.

\acknowledgements

This paper grew out of earlier work on local prediction for fields on regular
lattices with Rob Haslinger, Kristina Shalizi and Jacob Usinowicz; it is a
pleasure to acknowledge their support.  I have also benefited from discussions
with Jim Crutchfield, Dave Feldman, Christian Lindgren, Cris Moore, Mats
Nordahl, Mitch Porter and Karl Young.  Curtis Asplund pointed out typos in the
preprint.  I am grateful to Scott Page for the generous construction he puts on
the phrase ``complexity and diversity'', and to the organizers of DMCS 2003 for
the opportunity to participate.  Finally, I wish to thank Kara Kedi and
Kristina Shalizi for being perfect.

\bibliographystyle{crs}
\bibliography{locusts}

\end{document}